\documentclass[12pt,a4paper]{article}

\usepackage{amsmath}
\usepackage{theorem}
\usepackage{amsfonts}
\usepackage{amssymb}

\providecommand{\U}[1]{\protect\rule{.1in}{.1in}}

\newtheorem{theorem}{Theorem}

\newtheorem{lemma}[theorem]{Lemma}
\newtheorem{proposition}[theorem]{Proposition}
\newtheorem{remark}[theorem]{Remark}
\makeatletter
\def\blfootnote{\xdef\@thefnmark{}\@footnotetext}
\makeatother
\newenvironment{proof}[1][Proof]{\textit{#1.} }{\ \rule{0.5em}{0.5em}}

\begin{document}

\title{On the Location of the 1-particle Branch of the Spectrum of the
Disordered Stochastic Ising Model}
\author{M. Gianfelice \\
Dipartimento di Matematica \\
Universit\`{a} della Calabria \\
Ponte Pietro Bucci, cubo 30B, 87036 Arcavacata di Rende Italy \and M. Isopi 
\\
Dipartimento di Matematica \\
``Sapienza'' Universit\`{a} di Roma \\
P.le Aldo Moro 5, 00185 Roma Italy}
\maketitle

\begin{abstract}
We analyse the lower non trivial part of the spectrum of the generator of
the Glauber dynamics for a $d$-dimensional nearest neighbour Ising model
with a bounded random potential. We prove conjecture 1 in \cite{AMSZ}: for
sufficently large values of the temperature, the first band of the spectrum
of the generator of the process coincides with a closed non random segment
of the real line.
\end{abstract}

\blfootnote{
\par
\  \hspace*{0.2cm}\emph{Keywords} : Glauber dynamics, spectral gap, disordered
Ising model.}

\section{Introduction}

In \cite{AMSZ} the authors study the generator of the Glauber dynamics for a
one dimension Ising model with random bounded potential. They prove that, for
any realization of the potential and any value of the inverse temperature
$\beta>0,$ the spectrum of the generator is the union of disjoint closed
subsets of the real line ($k$-particle branches, $k\in\mathbb{N}^{+}$) and that,
with probability one with respect to the distribution of the potential, is a
non random set. In particular it is proved there that there exists a spectral
gap and thus the model exhibits exponential relaxation to equilibrium. As is
to be expected, and proved in \cite{AMSZ}, a relaxation rate which is valid
for every realization is the same as that of the non-disordered model with a
coupling constant that coincides with the maximum value of the coupling in the
disordered model. For the average over the disorder of the single spin
autocorrelation function, the speed of relaxation is somewhat larger as was
proved in \cite{Zh}.

Boundedness of the potential is essential for all these results of fast
convergence to equilibrium. In this case fairly detailed information on parts
of the spectrum of the generator is available (\cite{AMSZ}, \cite{Zh}). Also
in more than one dimension convergence slightly slower than exponential on
average can be proved at high temperature \cite{CMM}.

When the interactions are not bounded the situation is considerably different.
Even in one dimension there is no spectral gap (see \cite{Ze}) and relaxation
rate is subexponential (see \cite{SZ}).

In \cite{AMSZ} it is conjectured (conjecture 1, page 657) that results similar
to those proved there for one dimension should hold for $\beta$ small enough
in dimensions $d\geq2.$ It can be readily seen that for the proof, in one
dimension, of the results conjectured to be true in $d\geq2,$ the assumption
of ferromagnetic coupling is not needed. It is only used later to prove
exponential decay of eigenfunctions.

In this work we consider the Glauber dynamics for the $d$-dimensional nearest
neighbour Ising model, with a bounded random potential having absolutely
continuous distribution with respect to the Lebesgue measure and prove that
conjecture 1 in \cite{AMSZ} is true.

That is, there exists a constant $C,$ depending on the distribution of the
potential and on the lattice dimension $d$, such that, at high temperature,
the first branch of the spectrum of the generator of the process, at first
order in $\beta,$ coincides, for almost every realization of the potential,
with the segment
\[
\left[  1-C\beta,1+C\beta\right]
\]
(for a more precise statement see Theorem \ref{ct}). In particular this
implies that, at first order in $\beta,$ the spectral gap is larger than
$1-C\beta.$

We remark that at lower temperatures, but still in the uniqueness region,
relaxation is strictly slower than exponential for almost every realization of
the potential (see Theorem 3.3 of \cite{CMM}).

\section{Notations and results}

Consider the lattice $\mathbb{Z}^{d}$ and the set of bonds of the lattice
$\mathbb{B}_{d}:=\left\{  \left\{  x,y\right\}  \subset\mathbb{Z}%
^{d}:\left\vert x-y\right\vert =1\right\}  .$ We introduce a collection of
i.i.d random variables indexed by $\mathbb{B}_{d}.$ On each bond of the
lattice we define a random variable
\[
\omega_{b}\in\left[  J^{-},J^{+}\right]  ,\ b\in\mathbb{B}_{d}\ ,
\]
whose probability distribution is absolutely continuous with respect to the
Lebesgue measure. The random field $\omega$ is a function on the probability
space $\left(  \Omega,\mathcal{F},\mathbb{P}\right),\linebreak \Omega=\left[
J^{-},J^{+}\right]  ^{\mathbb{B}_{d}},$ and is ergodic w.r.t. the the group of
automorphisms on $\Omega$ generated by the lattice shift $\{\theta_{z}%
\}_{z\in\mathbb{Z}^{d}}$
\[
\Omega\ni\omega\longmapsto\omega\in\Omega:\left(  \theta_{z}\omega\right)
_{b}=\omega_{b-z}\quad z\in\mathbb{Z}^{d},b\in\mathbb{B}_{d}\,,
\]
where, $\forall\Lambda\subset\mathbb{Z}^{d},z\in\mathbb{Z}^{d},\,\Lambda
-z:=\{y\in\mathbb{Z}^{d}:y=x-z\,,\,x\in\Lambda\}\subset\mathbb{Z}^{d}.$

We now consider an Ising spin system in $\mathbb{Z}^{d}$. Denoting by
$\mathcal{S}$ the spin configuration space $\{-1,+1\}^{\mathbb{Z}^{d}}$and by
$\sigma$ the spin configuration, let $\left\{  \tau_{z}\right\}
_{z\in\mathbb{Z}^{d}}$ be the group of automorphisms of $\mathcal{S}$,
generated by the lattice translations
\[
\mathcal{S}\ni\sigma\longmapsto\tau_{z}\sigma\in\mathcal{S}:\left(  \tau
_{z}\sigma\right)  _{x}=\sigma_{x-z}\quad x,z\in\mathbb{Z}^{d}%
\]
and $j$ be the involution of $\mathcal{S}$ given by
\[
\mathcal{S}\ni\sigma\longmapsto j\left(  \sigma\right)  =-\sigma\in
\mathcal{S}\,.
\]

Let $\Lambda$ be a finite subset of the lattice. The Hamiltonian of the models
studied throughout this paper is
\begin{equation}
H_{\Lambda}^{\omega}\left(  \eta|\xi_{\partial\Lambda}\right)  =-\sum
_{x,y\in\Lambda\,:\,\left|  x-y\right|  =1}\frac{1}{2}\eta_{x}\omega_{x,y}%
\eta_{y}+\sum_{x\in\Lambda,\,y\in\Lambda^{c}\,:\,\left|  x-y\right|  =1}%
\eta_{x}\omega_{x,y}\xi_{y}\,, \label{defH}%
\end{equation}
where $\eta\in\mathcal{S}_{\Lambda}=\{-1,+1\}^{\Lambda}$ and
$\xi_{\partial\Lambda}:=\left(  \xi_{i}\right)  _{i\in\partial\Lambda}$
is a fixed boundary condition.

For any $\beta>0$ and any realization $\omega$ of the potential, let
$\mathcal{G}\left(  \beta,\omega\right)  $ be the set of Gibbs states of the
system specified by
\begin{align*}
\mu_{\Lambda}^{\beta,\omega}\left(  d\eta|\sigma_{\partial\Lambda}\right)   &
:=\frac{e^{-\beta H_{\Lambda}^{\omega}\left(  \eta|\sigma_{\partial\Lambda
}\right)  }}{Z_{\Lambda}^{\left(  d\right)  }\left(  \beta,\omega
|\sigma_{\partial\Lambda}\right)  }\mu_{\Lambda}\left(  d\eta\right)
\quad\Lambda\subset\subset\mathbb{Z}^{d}\\
Z_{\Lambda}^{\left(  d\right)  }\left(  \beta,\omega|\sigma_{\partial\Lambda
}\right)   &  :=\mu_{\Lambda}\left(  e^{-\beta H_{\Lambda}^{\omega}\left(
\eta|\sigma_{\partial\Lambda}\right)  }\right)  \,.
\end{align*}
We remark that, for a fixed boundary condition $\xi_{\partial\Lambda},$ the
conditional probability measure $\mu_{\Lambda}^{\beta,\omega}\left(  d\eta
|\xi_{\partial\Lambda}\right)  $ coincides with the one associated with the
formal Hamiltonian
\begin{equation}
H^{\omega}\left(  \sigma\right)  :=-\sum_{x,y\in\mathbb{Z}^{d}\,:\,\left\vert
x-y\right\vert =1}\frac{1}{2}\sigma_{x}\omega_{x,y}\sigma_{y}\,. \label{fH}%
\end{equation}
The Glauber processes studied in this paper are defined through the generator
\begin{equation}
L\left(  \beta,\omega\right)  f(\sigma):=\sum_{x\in\mathbb{Z}^{d}}w_{x}%
^{\beta,\omega}\left(  \sigma\right)  \left[  f\left(  \sigma\right)
-f\left(  \sigma^{x}\right)  \right]  \,, \label{defgen}%
\end{equation}
where the rates $w_{x}^{\beta,\omega}$ are chosen so that the process is
reversible w.r.t. $\mathcal{G}\left(  \beta,\omega\right)  $ and where
$\sigma^{x}$ represents the configuration in $\mathcal{S}$ such that
\[
\sigma_{y}^{x}=\left\{
\begin{array}
[c]{c}%
\sigma_{y}\quad y\neq x\\
-\sigma_{y}\quad y=x
\end{array}
\right.  \quad y\in\mathbb{Z}^{d}%
\]
and $f$ is a cylindrical function in $L^{2}\left(  \mathcal{S},\mu
^{\beta,\omega}\right)  :=\mathcal{L}\left(  \beta,\omega\right)  .$

We will always consider the generator $L$ a positive operator, so that \linebreak
$S\left(  t\right)  =\exp\left[  -tL\right]  $ will represent the associated semigroup.

In the following, with a little abuse of notation, we will use the same symbol
for the operator (\ref{defgen}) and for its closure in $\mathcal{L}\left(
\beta,\omega\right)  $ which, by reversibility of the Gibbs measure, is also
selfadjoint on $\mathcal{L}\left(  \beta,\omega\right)  $.

Let us define $J:=\left\vert J^{-}\right\vert \vee\left\vert J^{+}\right\vert
$ and, by (\ref{fH}), $\forall x\in\mathbb{Z}^{d},\,\omega\in\Omega$
\begin{equation}
\Delta_{x}H^{\omega}\left(  \sigma\right)  :=H^{\omega}\left(  \sigma\right)
-H^{\omega}\left(  \sigma^{x}\right)  =-\sigma_{x}\sum_{y:\,\left\vert
x-y\right\vert =1}\omega_{x,y}\sigma_{y}\,. \label{dHs}%
\end{equation}
Then
\begin{equation}
-4dJ\leq\left\vert \Delta_{x}H^{\omega}\left(  \sigma\right)  \right\vert
\leq4dJ\,. \label{maxdH}%
\end{equation}
{}From now on we are only interested in differences such as those in formula
(\ref{dHs}), which, as long as $x\in\Lambda,$ is the same regardless of
whether we use (\ref{defH}) or (\ref{fH}). So for simplicity we will be using
(\ref{fH}).

In the following we will restrict ourselves to the choice of transition rates
from $\sigma$ to $\sigma^{x}$ of the form
\begin{equation}
w_{x}^{\beta,\omega}\left(  \sigma\right)  =\psi\left(  \beta\Delta
_{x}H^{\omega}\left(  \sigma\right)  \right)  \,, \label{defw}%
\end{equation}
where $\psi$ is a monotone function, so that
\begin{equation}
\psi\left(  -\beta4dJ\right)  \wedge\psi\left(  \beta4dJ\right)  \leq
w_{x}^{\beta,\omega}\left(  \sigma\right)  \leq\psi\left(  -\beta4dJ\right)
\vee\psi\left(  \beta4dJ\right)  \,. \label{stimw}%
\end{equation}
In particular, we will work out the details for the case of the \emph{heat
bath dynamics} as was done in \cite{AMSZ}
\begin{equation}
w_{hb,x}^{\beta,\omega}\left(  \sigma\right)  =\psi_{hb}\left(  \beta
\Delta_{x}H^{\omega}\left(  \sigma\right)  \right)  =\frac{1}{1+e^{-\beta
\Delta_{x}H^{\omega}\left(  \sigma\right)  }}\ . \label{hbr}%
\end{equation}
Our analysis can be applied to any Glauber process with transition rates of
the kind given in (\ref{defw}).

The results contained in this paper are:

\begin{theorem}
\label{it}There exists a value $\beta_{d}^{-1}\left(  J\right)  $ of the
temperature such that, for any $\beta\in\left[  0,\beta_{d}\left(  J\right)
\right)  $ and any realization of the potential $\omega$, the first non
trivial branch of the spectrum of the generator of the heat bath dynamics,
$\sigma_{\beta}^{\left(  1\right)  }$, is contained in the interval $\left[
g_{d}^{-}\left(  \beta\right)  ,g_{d}^{+}\left(  \beta\right)  \right]  $
where $g_{d}^{-}\left(  \beta\right)  ,g_{d}^{+}\left(  \beta\right)  $ are
analytic functions of $\beta$ such that
\[
g_{d}^{\pm}\left(  \beta\right)  =1\pm2dJ\beta+o\left(  \beta\right)  \ .
\]

\end{theorem}

For a definition of $\sigma_{\beta}^{\left(  1\right)  }$ and a discussion of
its relevance see Corollary 1 of \cite{AMSZ} and Theorem 2.3 of \cite{M}.

\begin{theorem}
\label{mt}There exists a value $\beta_{d}^{\left(  1\right)  }$ of $\beta$
such that, for every $\beta\in\left[  0,\beta_{d}^{\left(  1\right)  }\right)
$ and almost every realization of the potential $\omega,$ the first non
trivial branch of the spectrum of the generator $\sigma_{\beta}^{\left(
1\right)  }$ satisfies
\[
\left[  1-f_{d}^{-}\left(  \beta\right)  ,1+f_{d}^{+}\left(  \beta\right)
\right]  \subseteq\sigma_{\beta}^{\left(  1\right)  }\,,
\]
where $f_{d}^{-}\left(  \beta\right)  ,\,f_{d}^{+}\left(  \beta\right)  $ are
analytic functions of $\beta$ such that
\[
f_{d}^{\pm}\left(  \beta\right)  =\pm2dJ\beta+o\left(  \beta\right)  \ .
\]

\end{theorem}

\begin{remark}
The analyticity of the functions introduced in the above two theorems, does
not hold only for the heat bath dynamics, but is guaranteed for any dynamics
where $\psi$ is an analytic function. If this is not the case, the statement
about analyticity must be dropped from the above theorems.
\end{remark}

\begin{theorem}
\label{ct}There exists a value $\beta_{d}^{\ast}\left(  J\right)  \leq
\beta_{d}^{\left(  1\right)  }\wedge\beta_{d}\left(  J\right)  $ of $\beta$
such that, for every $\beta\in\left[  0,\beta_{d}^{\ast}\left(  J\right)
\right)  $ and almost every realization of the potential $\omega$, the first
non trivial branch of the spectrum of the generator of the process
$\sigma_{\beta}^{\left(  1\right)  }$ is a non random set which coincides with
the closed subset of the real line $\left[  1-h_{d}^{-}(\beta),1+h_{d}%
^{+}(\beta)\right]  ,$ where $h_{d}^{\pm}\left(  \beta\right)  =\pm
2dJ\beta+o\left(  \beta\right)  .$
\end{theorem}

The proofs of these theorems rely in part on the approach of \cite{AMSZ} and
\cite{M} and in part on the lattice gas representation of the system, which we
will introduce in the next subsection. More precisely, we will restate the
dynamics with rates of the kind (\ref{defw}) in terms of a birth and death
process on the set of subsets of the lattice $\mathcal{P},$ which is naturally
isomorphic to $\mathcal{S},$ and make use of the setup given in \cite{GI1} and
\cite{GI2}.

\subsection{Lattice gas setting}

In \cite{GI1,GI2} we analysed the stochastic dynamics of a system with a
ferromagnetic potential constant on $\mathbb{B}_{d}$, confined in a finite
subset $\Lambda$ of the lattice and subject to free or periodic boundary
condition. Making use of a formalism borrowed from quantum mechanics, we were
able to represent the restriction of (\ref{defgen}) to $\mathcal{S}_{\Lambda
},$ in terms of a selfadjoint operator on $\mathcal{H}_{\Lambda}:=l^{2}\left(
\mathcal{P}_{\Lambda}\right)  $ which we showed to be unitarily equivalent
to a generator of birth and death process on $\mathcal{P}_{\Lambda}.$
Here we will follow the same approach.

We consider the Hilbert space of complex square summable function on the
single site configuration space with respect to the symmetric Bernoulli
measure. Namely, $\forall x\in{\mathbb{Z}}^{d},$%
\begin{align*}
\mathcal{H}_{x} &  :=span\left\{  \left\vert \emptyset\right\rangle
_{x},\left\vert \text{x}\right\rangle _{x}\right\}  \cong{\mathbb{C}}^{2}\\
\left\vert \emptyset\right\rangle _{x} &  \equiv\left(
\begin{array}
[c]{c}%
1\\
0
\end{array}
\right)  _{x}\quad\left\vert \text{x}\right\rangle _{x}\equiv\left(
\begin{array}
[c]{c}%
0\\
1
\end{array}
\right)  _{x}%
\end{align*}
$\mathcal{U}_{x}=M\left(  2,{\mathbb{C}}\right)  $ is the algebra of bounded
operators on $\mathcal{H}_{x}$\footnote[1]{Here, we think of $\mathcal{H}_{x}$
as spanned by two (orthonormal) vectors labelled by the \textquotedblright
empty site\textquotedblright\ and the \textquotedblright full
site\textquotedblright\ configurations. Consequently any operator acting on
the configuration space is lifted to a linear operator acting on
$\mathcal{H}_{x}$ and a probability density on the configuration space becomes
a convex combination of the projectors on the subspaces spanned by the basis
vectors of $\mathcal{H}_{x}$.}. Let us define the \textit{spin} operator
\[
\mathbf{s}_{x}\in\mathcal{U}_{x}:\mathbf{s}_{x}\left\{
\begin{array}
[c]{c}%
\left\vert \emptyset\right\rangle _{x}=\left\vert \text{x}\right\rangle _{x}\\
\left\vert \text{x}\right\rangle _{x}=\left\vert \emptyset\right\rangle _{x}%
\end{array}
\right.
\]
equivalent to the Pauli matrix $\sigma^{\left(  1\right)  }$%
\[
\sigma^{\left(  1\right)  }\equiv\left(
\begin{array}
[c]{cc}%
0 & 1\\
1 & 0
\end{array}
\right)
\]
and the \textit{spin flip} operator
\[
\mathbf{f}_{x}\in\mathcal{U}_{x}:\mathbf{f}_{x}\left\{
\begin{array}
[c]{c}%
\left\vert \emptyset\right\rangle _{x}=\left\vert \emptyset\right\rangle
_{x}\\
\left\vert \text{x}\right\rangle _{x}=-\left\vert \text{x}\right\rangle _{x}%
\end{array}
\right.
\]
equivalent to the Pauli matrix $\sigma^{\left(  3\right)  }$
\[
\sigma^{\left(  3\right)  }\equiv\left(
\begin{array}
[c]{cc}%
1 & 0\\
0 & -1
\end{array}
\right)  \ .
\]
Let $\Lambda$ be any finite subset of the $\mathbb{Z}^{d}$ lattice. Then we
have
\begin{align*}
\left\vert \alpha\right\rangle _{\Lambda} &  =\bigotimes_{x\in\alpha
}\left\vert \text{x}\right\rangle _{x}\bigotimes_{x\in\Lambda\backslash\alpha
}\left\vert \emptyset\right\rangle _{x}\\
\mathcal{H}_{\Lambda} &  =span\left\{  \left\vert \alpha\right\rangle
_{\Lambda}:\alpha\subseteq\Lambda\right\}
\end{align*}
Moreover $\mathcal{U}_{\Lambda}=M\left(  2^{\left\vert \Lambda\right\vert
},\mathbb{C}\right)  $ and $\mathcal{C}_{\Lambda}$ is the algebra of
polynomials in $\mathbf{s}_{\alpha}$ ($\mathbf{f}_{\alpha}$) $\forall
\alpha\subset\Lambda.$ Then
\begin{align*}
\mathbf{s}_{\alpha} &  =\bigotimes_{x\in\alpha}\mathbf{s}_{x}\bigotimes
_{x\in\Lambda\backslash\alpha}\mathbf{I}_{x}\ ,\\
\mathbf{f}_{\alpha} &  =\bigotimes_{x\in\alpha}\mathbf{f}_{x}\bigotimes
_{x\in\Lambda\backslash\alpha}\mathbf{I}_{x}\ ,\\
\mathbf{s}_{\emptyset} &  =\mathbf{f}_{\emptyset}=\mathbf{I}_{\Lambda}\ .
\end{align*}

Now, the generator of any Glauber process on the lattice, which in this
representation we denote by $\breve{L},$ can be written in terms of the
operators defined above and its generic matrix element becomes
\begin{equation}
\left(  \breve{L}\delta_{\alpha}\right)  _{\eta}=\sum_{x\in\mathbb{Z}^{d}%
}\left[  w\left(  \alpha,\alpha\triangle\{x\}\right)  \delta_{\eta,\alpha
}-w\left(  \alpha\triangle\{x\},\alpha\right)  \delta_{\eta,\alpha
\triangle\{x\}}\right]  \,,\label{meL}%
\end{equation}
where $\forall\alpha,\gamma\in\Lambda,\,\alpha\triangle\gamma=\left(
\alpha\cup\gamma\right)  \backslash\left(  \alpha\cap\gamma\right)  $ and,
with an abuse of notation, we indicate by $w\left(  \alpha,\alpha
\triangle\{x\}\right)  $ the transition rate from the state $\alpha$ to the
state\ $\alpha\triangle\{x\}.$

Since this form of the generator may seem unusual at first glance, here we
prove its equivalence to the classical form of generators of birth and death
processes on $\mathcal{P}.$

\subsubsection{Some remarks on birth and death processes for lattice gases}

We denote by $\mathbb{L}\left(  \mathcal{P}\right)  $ the linear space of
cylinder functions on $\mathcal{P}$ generated by linear combinations of
indicator functions of finite subsets of the lattice
\begin{align*}
\mathbb{L}\left(  \mathcal{P}\right)   &  \ni\varphi=\sum_{\alpha
\subset\mathbb{Z}^{d}\,:\,\left\vert \alpha\right\vert <\infty}\varphi
_{\alpha}\delta_{\alpha}\\
\forall\alpha &  \in\mathcal{P},\quad\mathcal{P}\ni\eta\longmapsto
\delta_{\alpha}\left(  \eta\right)  =\delta_{a,\eta}\in\left\{  0,1\right\}
\,,
\end{align*}
where the coefficients $\varphi_{\alpha}$ are real numbers.

Usually, see for example \cite{P}, the action of the generator of a birth and
death process $L$ on $\mathbb{L}\left(  \mathcal{P}\right)  $ takes a form
which can be expressed in either of the following two representations:
\begin{align}
L^{(-)}\varphi  :=\sum_{x\in\alpha}[ & w(\alpha\backslash\left\{  x\right\}  ,\alpha)  \left(  \varphi
_{\alpha\backslash\left\{  x\right\}  }-\varphi_{\alpha}\right)
\delta_{\alpha\backslash\{x\}} \label{Lbda-} \\
& +w(  \alpha,\alpha\backslash\left\{
x\right\})  \left( \varphi_{\alpha}-\varphi_{\alpha\backslash\left\{
x\right\}  }\right)  \delta_{\alpha}]  \nonumber \\
L^{(+)}\varphi  :=\sum_{x\in\alpha^{c}} [ & w(\alpha,\alpha\cup\left\{  x\right\})
 \left( \varphi_{\alpha}-\varphi_{\alpha\cup\left\{  x\right\}  }\right)  \delta_{\alpha}
\label{Lbda+} \\
 & +w(\alpha\cup\left\{  x\right\}  ,\alpha)  \left(  \varphi_{\alpha
\cup\left\{  x\right\}  }-\varphi_{\alpha}\right)  \delta_{\alpha\cup
\{x\}}] \nonumber \,.
\end{align}
Let $\mathcal{P}_{0}$ be the collection of finite and cofinite subsets of the
lattice. These expressions for $\left(  L\varphi\right)  _{\alpha}$ are
mutually equivalent and equivalent to
\[
\left(  \breve{L}\varphi\right)  _{\alpha}=\sum_{x\in\alpha}w\left(
\alpha,\alpha\backslash\left\{  x\right\}  \right)  \left(  \varphi_{\alpha
}-\varphi_{\alpha\backslash\left\{  x\right\}  }\right)  +\sum_{x\in\alpha
^{c}}w\left(  \alpha,\alpha\cup\left\{  x\right\}  \right)  \left(
\varphi_{\alpha}-\varphi_{\alpha\cup\left\{  x\right\}  }\right)  \ ,
\]
which can be derived from (\ref{meL}) (see (\ref{Lda}-\ref{Lbda}) below).
In fact, given the involution of $\mathcal{P}_{0}$
\begin{equation}
\mathcal{P}_{0}\ni\alpha\longmapsto\alpha^{c}=\mathbb{Z}^{d}\backslash
\alpha\in\mathcal{P}_{0}\quad\alpha\subset\mathbb{Z}^{d}\,,\label{definv}%
\end{equation}
we can define the family of operators $\left\{  \iota_{\Lambda}\right\}
_{\Lambda\in\mathcal{P}\,:\,\left\vert \Lambda\right\vert <\infty}$ on
$\mathbb{L}\left(  \mathcal{P}\right)  ,$ such that
\begin{align*}
\mathbb{L}\left(  \mathcal{P}\right)   &  \ni\varphi\longmapsto\phi
=\iota_{\Lambda}\varphi\in\mathbb{L}\left(  \mathcal{P}\right)  \\
\iota_{\Lambda}\delta_{\alpha} &  =\delta_{\alpha\triangle\Lambda}\quad
\alpha\in\mathcal{P}:\left\vert \alpha\right\vert <\infty\\
\iota_{\Lambda}\varphi &  =\sum_{\alpha\in\mathcal{P}\,:\,\left\vert
\alpha\right\vert <\infty}\varphi_{\alpha}\delta_{\alpha\triangle\Lambda}%
=\sum_{\alpha\in\mathcal{P}\,:\,\left\vert \alpha\right\vert <\infty}%
\varphi_{\alpha\triangle\Lambda}\delta_{\alpha}\,,\quad\varphi_{\alpha
\triangle\Lambda}=\varphi_{\left(  \alpha^{c}\cap\Lambda\right)  \cup\left(
\alpha\cap\Lambda^{c}\right)  }\\
\iota_{\Lambda}\left(  \iota_{\Lambda}\varphi\right)   &  =\varphi\quad
\varphi\in\mathbb{L}\left(  \mathcal{P}\right)  ,\,\Lambda\in\mathcal{P}%
:\left\vert \Lambda\right\vert <\infty\,.
\end{align*}
Defining $B$ to be the generator of a pure birth process with rates
\[
w\left(  \alpha\backslash\left\{  x\right\}  ,\alpha\right)  \mathbf{1}%
_{\alpha}\left(  x\right)  +w\left(  \alpha,\alpha\cup\left\{  x\right\}
\right)  \left(  1-\mathbf{1}_{\alpha}\left(  x\right)  \right)
\]
and $D$ the generator of a pure death process with rates
\[
w\left(  \alpha,\alpha\backslash\left\{  x\right\}  \right)  \mathbf{1}%
_{\alpha}\left(  x\right)  +w\left(  \alpha\cup\left\{  x\right\}
,\alpha\right)  \left(  1-\mathbf{1}_{\alpha}\left(  x\right)  \right)  \,,
\]
we may rewrite (\ref{Lbda-}) and (\ref{Lbda+}) in the form
\[
\left(  L^{\left(  \pm\right)  }\varphi\right)  _{\alpha}=\left(  B^{\left(
\pm\right)  }\varphi\right)  _{\alpha}+\left(  D^{\left(  \pm\right)  }%
\varphi\right)  _{\alpha}\,,
\]
where the definition of $B^{\left(  \pm\right)  }$ and $D^{\left(  \pm\right)
}$ is readily understood. Since
\begin{align}
w\left(  \alpha\backslash\left\{  x\right\}  ,\alpha\right)   &  =w\left(
\alpha^{c}\cup\left\{  x\right\}  ,\alpha^{c}\right)  \label{bvsd}\\
w\left(  \alpha,\alpha\backslash\left\{  x\right\}  \right)   &  =w\left(
\alpha^{c},\alpha^{c}\cup\left\{  x\right\}  \right)  \,,\nonumber
\end{align}
considering for example (\ref{Lbda-}), for any finite $\Lambda\subset
\mathbb{Z}^{d}$ we have
\begin{align*}
\left(  \iota_{\Lambda}B^{\left(  -\right)  }\iota_{\Lambda}\varphi\right)
_{\alpha} &  =\sum_{x\in\alpha\triangle\Lambda}w\left(  \left(  \alpha
\triangle\Lambda\right)  \backslash\left\{  x\right\}  ,\alpha\triangle
\Lambda\right)  \left(  \left(  \iota_{\Lambda}\varphi\right)  _{\left(
\left(  \alpha\triangle\Lambda\right)  \backslash\left\{  x\right\}  \right)
}-\left(  \iota_{\Lambda}\varphi\right)  _{\left(  \alpha\triangle
\Lambda\right)  }\right)  \\
&  =\sum_{x\in\alpha\triangle\Lambda}w\left(  \left(  \alpha\triangle
\Lambda\right)  \backslash\left\{  x\right\}  ,\alpha\triangle\Lambda\right)
\left(  \varphi_{\left(  \left(  \alpha\triangle\Lambda\right)  \backslash
\left\{  x\right\}  \right)  \triangle\Lambda}-\varphi_{\alpha}\right)
\end{align*}
and choosing $\Lambda\supset\alpha,$ by (\ref{bvsd}) we get
\begin{align*}
\left(  \iota_{\Lambda}B^{\left(  -\right)  }\iota_{\Lambda}\varphi\right)
_{\alpha} &  =\sum_{x\in\alpha^{c}\cap\Lambda}w\left(  \left(  \alpha^{c}%
\cap\Lambda\right)  \backslash\left\{  x\right\}  ,\alpha^{c}\cap
\Lambda\right)  \left(  \varphi_{\left(  \left(  \alpha^{c}\cap\Lambda\right)
\backslash\left\{  x\right\}  \right)  ^{c}\cap\Lambda}-\varphi_{\alpha
}\right)  \\
&  =\sum_{x\in\alpha^{c}\cap\Lambda}w\left(  \left(  \alpha\cup\left\{
x\right\}  \right)  ^{c}\cap\Lambda,\alpha^{c}\cap\Lambda\right)  \left(
\varphi_{\alpha\cup\left\{  x\right\}  }-\varphi_{\alpha}\right)  \\
&  =\sum_{x\in\alpha^{c}\cap\Lambda}w\left(  \alpha\cup\left\{  x\right\}
\cup\Lambda^{c},\alpha\cup\Lambda^{c}\right)  \left(  \varphi_{\alpha
\cup\left\{  x\right\}  }-\varphi_{\alpha}\right)  \\
&  =\sum_{x\in\left(  \alpha\cup\Lambda^{c}\right)  ^{c}}w\left(  \left(
\alpha\cup\Lambda^{c}\right)  \cup\left\{  x\right\}  ,\alpha\cup\Lambda
^{c}\right)  \left(  \varphi_{\alpha\cup\left\{  x\right\}  }-\varphi_{\alpha
}\right)  \,.
\end{align*}
We now assume the system to be confined in a box $\Lambda$ with boundary
conditions $\eta.$ Let $\mathcal{P}_{\Lambda}$ be the set of the subsets of
$\Lambda.$ We can inject $\mathbb{L}\left(  \mathcal{P}_{\Lambda}\right)  ,$
the vector space generated by linear combinations of $\delta_{\alpha}%
,\ \alpha\subseteq\Lambda,$ in $\mathbb{L}\left(  \mathcal{P}\right)  $ and
consider a naturally defined $\iota_{\Lambda}^{\eta}.$
\begin{align}
\mathbb{L}\left(  \mathcal{P}_{\Lambda}\right)   &  \ni\varphi\longmapsto
\phi^{\eta}=\iota_{\Lambda}^{\eta}\varphi=\iota_{\Lambda}\left(  \varphi
\delta_{\eta}\right)  \in\mathbb{L}\left(  \mathcal{P}\right)  \label{defi}\\
\iota_{\Lambda}^{\eta}\delta_{\alpha} &  =\delta_{\left(  \alpha\cup
\eta\right)  \triangle\Lambda}=\delta_{\Lambda\backslash\alpha\cup\eta}%
\quad\alpha\subseteq\Lambda\,.
\end{align}

Independently of the choice of the boundary conditions $\eta,$ $\forall
\alpha\subseteq\Lambda$
\begin{eqnarray*}
\left(  \iota_{\Lambda}^{\eta}B_{\Lambda,\eta}^{\left(  -\right)  }%
\iota_{\Lambda}^{\eta}\varphi\right)  _{\alpha} & = & \sum_{x\in\Lambda
\backslash\alpha}w\left(  \alpha\cup\left\{  x\right\}  ,\alpha\right)
\left(  \varphi_{\alpha\cup\left\{  x\right\}  }-\varphi_{\alpha}\right) \\
& = & \left(  D_{\Lambda,\eta}^{\left(  +\right)  }\varphi\right)  _{\alpha}%
\quad\eta\in\mathcal{P}_{\Lambda^{c}},\,\left\vert \eta\right\vert <\infty\,,
\end{eqnarray*}
where $B_{\Lambda,\eta}^{\left(  \pm\right)  }$ and $D_{\Lambda,\eta}^{\left(
\pm\right)  }$ denote the natural restrictions of $B^{\left(  \pm\right)  }$
and $D^{\left(  \pm\right)  },\,$to $\mathbb{L}\left(  \mathcal{P}_{\Lambda
}\right)  .$

To keep notation simple, from now on we will omit to indicate the boundary
conditions where there is no risk of ambiguity.

Now, for any realization of $\omega\in\Omega,\ \beta\geq0$ and $\alpha
\in\mathcal{P},$ denoting by $\bar{L}\left(  \beta,\omega\right)  $ the
generator of the process given in (\ref{defgen}) in this representation, from
(\ref{meL}) we get
\begin{align}
\bar{L}\left(  \beta,\omega\right)  \delta_{\alpha} &  =\sum_{\eta
\in\mathcal{P}\,:\,\left\vert \eta\right\vert <\infty}\left(  \bar{L}\left(
\beta,\omega\right)  \delta_{\alpha}\right)  _{\eta}\delta_{\eta}\label{Lda}\\
&  =\sum_{x\in\mathbb{Z}^{d}}\left[  w^{\beta,\omega}\left(  \alpha
,\alpha\triangle\{x\}\right)  \delta_{\alpha}-w^{\beta,\omega}\left(
\alpha\triangle\{x\},\alpha\right)  \delta_{\alpha\triangle\{x\}}\right]
\,.\nonumber
\end{align}
Then, $\forall\varphi\in\mathbb{L}\left(  \mathcal{P}\right)  ,$ we have
\begin{align}
\bar{L}\left(  \beta,\omega\right)  \varphi &  =\sum_{\alpha\in\mathcal{P}%
\,:\,\left\vert \alpha\right\vert <\infty}\left(  \bar{L}\left(  \beta
,\omega\right)  \varphi\right)  _{\alpha}\delta_{\alpha}\label{Lbd}\\
\bar{L}\left(  \beta,\omega\right)  \varphi &  =\sum_{\alpha\in\mathcal{P}%
\,:\,\left\vert \alpha\right\vert <\infty}\sum_{x\in\mathbb{Z}^{d}}%
\varphi_{\alpha}\left[  w^{\beta,\omega}\left(  \alpha,\alpha\triangle
\{x\}\right)  \delta_{\alpha}-w^{\beta,\omega}\left(  \alpha\triangle
\{x\},\alpha\right)  \delta_{\alpha\triangle\{x\}}\right]  \label{Llg}\\
&  =\sum_{\alpha\in\mathcal{P}\,:\,\left\vert \alpha\right\vert <\infty}%
\sum_{x\in\mathbb{Z}^{d}}w^{\beta,\omega}\left(  \alpha,\alpha\triangle
\{x\}\right)  \left(  \varphi_{a}-\varphi_{a\triangle\{x\}}\right)
\delta_{\alpha}\,.\nonumber
\end{align}
Notice that, for any $\omega\in\Omega,$ (\ref{Llg}) takes the form
\begin{align}
\left(  \bar{L}\left(  \beta,\omega\right)  \varphi\right)  _{\alpha}
:=&\sum_{x\in\alpha}w^{\beta,\omega}\left(  \alpha,\alpha\backslash\left\{
x\right\}  \right)  \left(  \varphi_{\alpha}-\varphi_{\alpha\backslash\left\{
x\right\}  }\right)  +\label{Lbda}\\
&  +\sum_{x\in\alpha^{c}}w^{\beta,\omega}\left(  \alpha,\alpha\cup\left\{
x\right\}  \right)  \left(  \varphi_{\alpha}-\varphi_{\alpha\cup\left\{
x\right\}  }\right)  \ ,\quad\alpha\in\mathcal{P}:\left\vert \alpha\right\vert
<\infty\ .\nonumber
\end{align}
Hence, since $\iota_{\Lambda}D_{\Lambda}^{\left(  +\right)  }\iota_{\Lambda
}=B_{\Lambda}^{\left(  -\right)  },$
\[
L_{\Lambda}^{\left(  -\right)  }=B_{\Lambda}^{\left(  -\right)  }+D_{\Lambda
}^{\left(  -\right)  }=\iota_{\Lambda}\left(  B_{\Lambda}^{\left(  +\right)
}+D_{\Lambda}^{\left(  +\right)  }\right)  \iota_{\Lambda}=\iota_{\Lambda
}L^{\left(  +\right)  }\iota_{\Lambda}%
\]
and, $\forall\omega\in\Omega,\ \beta\geq0,\,\alpha\subseteq\Lambda$ and
boundary condition $\eta,$ (\ref{Lbda}), takes the form
\begin{align}
\left(  \bar{L}_{\Lambda}\left(  \beta,\omega,\eta\right)  \varphi\right)
_{\alpha} &  =\left(  D_{\Lambda}^{\left(  -\right)  }\left(  \omega
,\beta,\eta\right)  \varphi\right)  _{\alpha}+\left(  B_{\Lambda}^{\left(
+\right)  }\left(  \omega,\beta,\eta\right)  \varphi\right)  _{\alpha
}\label{Lbdeq}\\
&  =\left(  D_{\Lambda}^{\left(  -\right)  }\left(  \omega,\beta,\eta\right)
\varphi\right)  _{\alpha}+\left(  \iota_{\Lambda}D_{\Lambda}^{\left(
-\right)  }\left(  \omega,\beta,\eta\right)  \iota_{\Lambda}\varphi\right)
_{\alpha}\,.\nonumber
\end{align}

It is worth to notice that, for any realization of the potential, $\bar
{L}_{\Lambda}\left(  \beta,\omega,\eta\right)  $ commutes with $\iota
_{\Lambda}.$

We remark that the equivalence between (\ref{Llg}) and the generator of
process defined in (\ref{defgen}) can be deduced comparing the associated
Dirichlet forms.

\section{Proof of the Theorems}

Replacing $\varphi$ by $\delta_{\eta}$ for a fixed $\eta\subseteq\Lambda$ in
(\ref{Lbda}), we get the generic matrix element of (\ref{Lbd}) and then of
(\ref{Llg}) as operators acting on $\mathcal{H}_{\Lambda}.$ We can then
transform (\ref{Llg}) into a selfadjoint operator $\tilde{L}_{\Lambda}%
^{s}\left(  \beta,\omega\right)  $ on $\mathcal{H}_{\Lambda}$ through the
unitary mapping from $\mathcal{H}_{\Lambda}\left(  \beta,\omega\right)
:=l^{2}\left(  \mathcal{P}_{\Lambda},\mu_{\Lambda}^{\beta,\omega}\right)  $
(which is isomorphic to the restriction of $\mathcal{L}\left(  \beta
,\omega\right)  $ to $\Lambda$) to $\mathcal{H}_{\Lambda}$\ given by the
multiplication of the elements of $\mathcal{H}_{\Lambda}\left(  \beta
,\omega\right)  $ by $\sqrt{\frac{\mu_{\Lambda}^{\beta,\omega}}{\mu_{\Lambda}%
}}.$

We will give a relative bound of the Dirichlet form of $\tilde{L}_{\Lambda
}^{s}\left(  \beta,\omega\right)  $ in terms of the Dirichlet form of the
generator of the independent process $\bar{L}_{\Lambda}$ and make use of
standard perturbation theory to give a lower bound for the spectral gap of
$\tilde{L}_{\Lambda}^{s}\left(  \beta,\omega\right)  ,\,g_{d}^{-,\Lambda
}\left(  \beta\right)  ,$ for small values of $\beta>0$ and for any $\omega
\in\Omega.$ These bounds will turn out to be independent of $\Lambda,$ which
implies in particular $g_{d}^{-,\Lambda}\left(  \beta\right)  =g_{d}%
^{-}\left(  \beta\right)  ,$ and therefore extend to the infinite volume
setting. We get $g_{d}^{+}\left(  \beta\right)  $ by applying the same
argument to the operator $\hat{L}_{\Lambda}^{s}\left(  \beta,\omega\right)
\geq\tilde{L}_{\Lambda}^{s}\left(  \beta,\omega\right)  $ on $\mathcal{H}%
_{\Lambda},$ which is also unitary equivalent to a generator of a Glauber
process for the Ising model reversible with respect to $\mu_{\Lambda}%
^{\beta,\omega}.$ The proof of Theorem 2 relies on two results. First a
theorem of Minlos (Theorem 2.2 of \cite{M}) which gives detailed information
on the first branch of the spectrum for constant realizations. Second on the
part 2) of Theorem 3 of \cite{AMSZ}, which proves that the first branch of the
spectrum for a constant realization is contained in the first branch of the
spectrum with random coupling.

Finally, since the family of operators and spaces $\left(  L\left(
\beta,\omega\right)  ,\mathcal{L}\left(  \beta,\omega\right)  \right)  $ is a
metrically transitive family with respect to lattice translations, from
general results of spectral theory for random operators (see \cite{PF} and
Remark 4 of \cite{AMSZ}), it will follow that the spectrum of $L\left(
\beta,\omega\right)  $ is non-random for $\mathbb{P}$-a.e. $\omega.$ This
remark, together with the two previous results, will then prove Theorem 3.

Let us consider the heat bath case. Given a finite portion of the lattice
$\Lambda$ and a realization of the potential $\omega,$ assuming for example
periodic boundary condition, the restriction of the generator of the process
given in (\ref{Llg}) to $\mathcal{P}_{\Lambda},$ takes the form (\ref{Lbd}),
(\ref{Lbda}), where
\[
w^{\beta,\omega}\left(  \alpha,\alpha\triangle\left\{  x\right\}  \right)
=\psi_{hb}\left(  \beta\Delta_{x}H_{\alpha}^{\omega}\right)  =\frac
{1}{1+e^{-\beta\Delta_{x}H_{\alpha}^{\omega}}}%
\]
and
\begin{align}
\Delta_{x}H_{\alpha}^{\omega}  &  =H_{\alpha}\left(  \omega\right)
-H_{\alpha\triangle\{x\}}\left(  \omega\right) \label{defgradH}\\
&  =\mathbf{1}_{\alpha}\left(  x\right)  [H_{\alpha}\left(  \omega\right)
-H_{\alpha\backslash\{x\}}\left(  \omega\right)  ]+\mathbf{1}_{\alpha^{c}%
}\left(  x\right)  [H_{\alpha}\left(  \omega\right)  -H_{\alpha\cup
\{x\}}\left(  \omega\right)  ]\nonumber
\end{align}
representing respectively (\ref{hbr}) and (\ref{dHs}) in the lattice gas
framework. Here
\[
H_{\alpha}(\omega)=\sum_{\ b\in\mathbb{B}_{d}}\omega_{b}-2\sum_{\ b\in
\partial\alpha}\omega_{b}.
\]
Although infinite, $\sum_{\ b\in\mathbb{B}_{d}}\omega_{b}$ is an harmless
constant since transition rates are functions only of $\Delta_{x}H_{\alpha
}^{\omega}$.

Following \cite{GI1}, since $\mathcal{H}_{\Lambda}\cong\bigoplus
_{n=0}^{\left\vert \Lambda\right\vert }\mathcal{H}_{\Lambda}^{\left(
n\right)  }$, with $\mathcal{H}_{\Lambda}^{\left(  0\right)  }\equiv
\mathbb{R}$ and
\[
\mathcal{H}_{\Lambda}^{\left(  n\right)  }:=\{\left\vert
\alpha\right\rangle \in\mathcal{H}_{\Lambda}:\left\vert \alpha\right\vert
=n\} \,,
\]
we denote by $U_{\Lambda}$ the unitary operator
\begin{align}
U_{\Lambda}  &  :\mathcal{H}_{\Lambda}\longrightarrow\mathcal{H}_{\Lambda
}\label{defU}\\
U_{\Lambda}\left\vert \alpha\right\rangle  &  =\frac{1}{2^{\frac{\left\vert
\Lambda\right\vert }{2}}}\sum_{\gamma\subseteq\Lambda}\left(  -1\right)
^{\left\vert \alpha\cap\gamma\right\vert }\left\vert \gamma\right\rangle
\quad\alpha\subseteq\Lambda\,,\nonumber
\end{align}
and by $E_{\Lambda},$ the representation of the involution $\iota_{\Lambda}$
introduced in (\ref{defi}) as an operator on $\mathcal{H}_{\Lambda},$ that is
\begin{align}
E_{\Lambda}  &  :\mathcal{H}_{\Lambda}\longrightarrow\mathcal{H}_{\Lambda
}\label{defE}\\
E_{\Lambda}\left\vert \alpha\right\rangle  &  =\left\vert \Lambda
\backslash\alpha\right\rangle \quad\alpha\subseteq\Lambda\,.\nonumber
\end{align}
Now, for any $\alpha\subseteq\Lambda,\ E_{\Lambda}\frac{\left\vert
\alpha\right\rangle \pm\left\vert \alpha^{c}\right\rangle }{\sqrt{2}}=\pm
\frac{\left\vert \alpha\right\rangle \pm\left\vert \alpha^{c}\right\rangle
}{\sqrt{2}},$ hence $\mathcal{H}_{\Lambda}=\mathcal{H}_{\Lambda}^{+}%
\oplus\mathcal{H}_{\Lambda}^{-},$ where
\[
\mathcal{H}_{\Lambda}^{\pm}:=span\{\frac{\left\vert \alpha\right\rangle
\pm\left\vert \alpha^{c}\right\rangle }{\sqrt{2}}:\alpha\subseteq\Lambda\}
\]
Moreover, setting
\begin{equation}
\bar{E}_{\Lambda}:=U_{\Lambda}E_{\Lambda}U_{\Lambda}:\mathcal{H}_{\Lambda
}\longrightarrow\mathcal{H}_{\Lambda}\,, \label{defEbar}%
\end{equation}
since
\begin{eqnarray}
\delta_{\alpha,\gamma} & = &\left\langle \alpha|\gamma\right\rangle =\left\langle
\alpha\right\vert U_{\Lambda}U_{\Lambda}\left\vert \gamma\right\rangle
=2^{-\left\vert \Lambda\right\vert }\sum_{\eta\subseteq\Lambda}\left(
-1\right)  ^{\left\vert \alpha\cap\eta\right\vert +\left\vert \gamma\cap
\eta\right\vert } \label{U2} \\
& =& 2^{-\left\vert \Lambda\right\vert }\sum_{\eta\subseteq
\Lambda}\left(  -1\right)  ^{\left\vert \left(  \alpha\triangle\gamma\right)
\cap\eta\right\vert } \nonumber \,, %
\end{eqnarray}
then, for any $\left\vert \alpha\right\rangle \in\mathcal{H}_{\Lambda},$
\[
\bar{E}_{\Lambda}\left\vert \alpha\right\rangle =\sum_{\gamma,\eta
\subseteq\Lambda}\frac{\left(  -1\right)  ^{\left\vert \alpha\cap
\gamma\right\vert +\left\vert \left(  \Lambda\backslash\gamma\right)  \cap
\eta\right\vert }}{2^{\left\vert \Lambda\right\vert }}\left\vert
\eta\right\rangle =\sum_{\gamma,\eta\subseteq\Lambda}\frac{\left(  -1\right)
^{\left\vert \alpha\cap\gamma\right\vert -\left\vert \gamma\cap\eta\right\vert
+\left\vert \eta\right\vert }}{2^{\left\vert \Lambda\right\vert }}\left\vert
\eta\right\rangle =\left(  -1\right)  ^{\left\vert \alpha\right\vert
}\left\vert \alpha\right\rangle \,,
\]
so that $\mathcal{H}_{\Lambda}$ can also be decomposed as  direct sum of
$\overline{\mathcal{H}}_{\Lambda}^{+}:=\bigoplus_{n\geq0\,:\,2n\in
\{0,..,\left\vert \Lambda\right\vert \}}\mathcal{H}_{\Lambda}^{\left(
2n\right)  }$ and $\overline{\mathcal{H}}_{\Lambda}^{-}:=\bigoplus
_{n\geq0\,:\,2n+1\in\{0,..,\left\vert \Lambda\right\vert \}}\mathcal{H}%
_{\Lambda}^{\left(  2n+1\right)  }.$ Clearly $U_{\Lambda}\mathcal{H}_{\Lambda
}^{\pm}=\overline{\mathcal{H}}_{\Lambda}^{\pm}.$ Furthermore, by
(\ref{U2}), $\mathcal{H}_{\Lambda}=\mathcal{H}_{\Lambda,+}%
\oplus\mathcal{H}_{\Lambda,-}\,,$ where $U_{\Lambda}\mathcal{H}_{\Lambda,\pm
}=\pm\mathcal{H}_{\Lambda,\pm}.$

If, for any $x\in\Lambda$, $\ell_{x}^{\Lambda},\,\ell_{x}^{\Lambda,\bot}$
denote the mutually orthogonal projectors on $\mathcal{H}_{\Lambda}$ such
that
\begin{equation}
\ell_{x}^{\Lambda,\bot}:=I_{\Lambda}-\ell_{x}^{\Lambda}\,;\,\ell_{x}^{\Lambda
}\left\vert \alpha\right\rangle =\mathbf{1}_{\alpha}\left(  x\right)
\left\vert \alpha\right\rangle \quad\alpha\subseteq\Lambda\,. \label{lid}%
\end{equation}
We have
\begin{gather}
\bar{\ell}_{x}^{\Lambda}=U_{\Lambda}\ell_{x}^{\Lambda}U_{\Lambda};\quad
\bar{\ell}_{x}^{\Lambda,\bot}=U_{\Lambda}\ell_{x}^{\Lambda,\bot}U_{\Lambda
},\label{lid1}\\
\ell_{x}^{\Lambda}=E_{\Lambda}\ell_{x}^{\Lambda,\bot}E_{\Lambda};\quad\ell
_{x}^{\Lambda,\bot}=E_{\Lambda}\ell_{x}^{\Lambda}E_{\Lambda},\nonumber\\
\lbrack E_{\Lambda},\bar{\ell}_{x}^{\Lambda}]=[\bar{E}_{\Lambda},\ell
_{x}^{\Lambda}]=0.\nonumber
\end{gather}
We also denote by
\[
\frac{e^{-\frac{\beta}{2}H_{\Lambda}^{\omega}}}{\left(  Z_{\Lambda}^{\left(
d\right)  }\left(  \beta,\omega\right)  \right)  ^{\frac{1}{2}}}%
:\mathcal{H}_{\Lambda}\left(  \beta,\omega\right)  \longrightarrow
\mathcal{H}_{\Lambda}%
\]
the matrix representation of the multiplication operator by $\sqrt{\frac
{\mu_{\Lambda}^{\beta,\omega}}{\mu_{\Lambda}}}.$

In \cite{GI1,GI2}, comparing Dirichlet forms, we also showed that $\bar
{L}_{\Lambda}^{s}\left(  \beta,\omega\right)  $ admits the representation
\[
\tilde{L}_{\Lambda}^{s}\left(  \beta,\omega\right)  =\sum_{x\in\Lambda}%
\tilde{L}_{x,\Lambda}^{s}\left(  \beta,\omega\right)  \,,
\]
whose matrix elements, by the definition of $\Delta_{x}H_{\alpha}^{\omega},$
are, for any two vectors $\left\vert \alpha\right\rangle ,\left\vert
\gamma\right\rangle $ of the basis of $\mathcal{H}_{\Lambda}\ $

\begin{align}
& \left\langle \gamma\right\vert \tilde{L}_{\Lambda}^{s}\left(  \beta
,\omega\right)  \left\vert \alpha\right\rangle  =\sum_{x\in\Lambda
}\left\langle \gamma\right\vert \tilde{L}_{x,\Lambda}^{s}\left(  \beta
,\omega\right)  \left\vert \alpha\right\rangle \label{LbdsI}\\
&\left\langle \gamma\right\vert \tilde{L}_{x,\Lambda}^{s}\left(  \beta
,\omega\right)  \left\vert \alpha\right\rangle  \nonumber\\  =& \left\langle
\gamma\right\vert \left\{  \mathbf{1}_{\alpha}\left(  x\right)  \frac{1}%
{\cosh\frac{\beta}{2}\Delta_{x}H_{\alpha}^{\omega}}\left[  \bar{\ell}%
_{x}^{\Lambda}+I_{\Lambda}\frac{e^{\frac{\beta}{2}\Delta_{x}H_{\alpha}%
^{\omega}}-1}{2}\right]  +\right. \nonumber\\
& \qquad\quad \left. +\mathbf{1}_{\alpha^{c}}\left(  x\right)  \frac{1}{\cosh\frac
{\beta}{2}\Delta_{x}H_{\alpha}^{\omega}}\left[  \bar{\ell}_{x}^{\Lambda
}+I_{\Lambda}\frac{e^{\frac{\beta}{2}\Delta_{x}H_{\alpha}^{\omega}}-1}%
{2}\right]  \right\}  \left\vert \alpha\right\rangle \nonumber\\
=& \left\langle \gamma\right\vert \left\{  \mathbf{1}_{\alpha}\left(
x\right)  \frac{1}{\cosh\frac{\beta}{2}\Delta_{x}H_{\alpha}^{\omega}}\left[
\frac{e^{\frac{\beta}{2}\Delta_{x}H_{\alpha}^{\omega}}+1}{2}\bar{\ell}%
_{x}^{\Lambda}+\frac{e^{\frac{\beta}{2}\Delta_{x}H_{\alpha}^{\omega}}-1}%
{2}\bar{\ell}_{x}^{\Lambda,\bot}\right]  +\right. \nonumber\\
& \qquad\quad \left.  +\mathbf{1}_{\alpha^{c}}\left(  x\right)  \frac{1}{\cosh\frac
{\beta}{2}\Delta_{x}H_{\alpha}^{\omega}}\left[  \frac{e^{\frac{\beta}{2}%
\Delta_{x}H_{\alpha}^{\omega}}+1}{2}\bar{\ell}_{x}^{\Lambda}+\frac
{e^{\frac{\beta}{2}\Delta_{x}H_{\alpha}^{\omega}}-1}{2}\bar{\ell}_{x}%
^{\Lambda,\bot}\right]  \right\}  \left\vert \alpha\right\rangle \,.\nonumber
\end{align}

\subsection{Proof of Theorem \ref{it}}

Let us set
\begin{align*}
L_{\Lambda}  &  :=\sum_{x\in\Lambda}\ell_{x}^{\Lambda}\,;\quad\bar{L}%
_{\Lambda}=U_{\Lambda}L_{\Lambda}U_{\Lambda}=\sum_{x\in\Lambda}\bar{\ell}%
_{x}^{\Lambda}\,,\\
L_{\Lambda}\left|  \alpha\right\rangle  &  =\left|  \alpha\right|  \left|
\alpha\right\rangle \quad\alpha\subseteq\Lambda\,.
\end{align*}

\begin{lemma}
For any $\left\vert u\right\rangle \in\mathcal{H}_{\Lambda}^{\pm},$
\begin{equation}
\left\langle u\right\vert \bar{L}_{\Lambda}\left\vert u\right\rangle
\leq2\left\langle u\right\vert L_{\Lambda}\left\vert u\right\rangle \ .
\label{Lbar<L}%
\end{equation}

\end{lemma}

\begin{proof}
We first notice that, for any $x\in\Lambda,\,a:\mathcal{P}_{\Lambda}%
\times\mathcal{P}_{\Lambda}\longrightarrow\mathbb{R},$
\[
\sum_{\alpha\subseteq\Lambda}a_{\alpha,\alpha\cup\{x\}}\mathbf{1}_{\alpha^{c}%
}\left(  x\right)  =\sum_{\alpha\subseteq\Lambda}\mathbf{1}_{\alpha}\left(
x\right)  a_{\alpha\backslash\{x\},\alpha}\,.
\]

Then, for any $\left\vert u\right\rangle \in\mathcal{H}_{\Lambda}^{\pm},$ we
get
\begin{align*}
\left\langle u\right\vert \bar{L}_{\Lambda}\left\vert u\right\rangle  &
=\sum_{x\in\Lambda}\sum_{\alpha\subseteq\Lambda}\frac{u_{\alpha}^{2}%
-u_{\alpha}u_{\alpha\triangle\{x\}}}{2}=\sum_{x\in\Lambda}\sum_{\alpha
\subseteq\Lambda}\left(  \frac{u_{\alpha}-u_{\alpha\triangle\{x\}}}{2}\right)
^{2}\\
&  \leq\sum_{x\in\Lambda}\sum_{\alpha\subseteq\Lambda}\frac{u_{\alpha}%
^{2}+u_{\alpha\triangle\{x\}}^{2}}{2}=\sum_{\alpha\subseteq\Lambda}\sum
_{x\in\alpha}\frac{u_{\alpha}^{2}+u_{\alpha\backslash\{x\}}^{2}}{2}%
+\sum_{\alpha\subseteq\Lambda}\sum_{x\in\alpha^{c}}\frac{u_{\alpha}%
^{2}+u_{\alpha\cup\{x\}}^{2}}{2}\\
&  =\sum_{\alpha\subseteq\Lambda}\sum_{x\in\alpha}u_{\alpha}^{2}+\sum
_{\alpha\subseteq\Lambda}\sum_{x\in\alpha^{c}}u_{\alpha}^{2}=\sum
_{\alpha\subseteq\Lambda}\sum_{x\in\alpha}u_{\alpha}^{2}+\sum_{\alpha
^{c}\subseteq\Lambda}\sum_{x\in\alpha^{c}}u_{\alpha^{c}}^{2}=2\left\langle
u\right\vert L_{\Lambda}\left\vert u\right\rangle \,.
\end{align*}

\end{proof}

\begin{remark}
{}From (\ref{defH}) it follows that, for any $\omega\in\Omega,\ H_{\alpha
}^{\omega}$ depends on $\alpha$ only through the subset of $\mathbb{B}%
_{\Lambda},$ $\partial\alpha:=\left\{  b\in\mathbb{B}_{\Lambda}:\left\vert
b\cap\alpha\right\vert =1\right\}  $ then, because $\partial\alpha
=\partial\alpha^{c}$, by (\ref{Lbdeq}), (\ref{defE}) and (\ref{lid1}), for any
realization of the potential $H_{\alpha}\left(  \omega\right)  =H_{\alpha^{c}%
}\left(  \omega\right)  .$ Hence, for any $\beta\geq0$ and $\omega\in
\Omega,\ \tilde{L}_{\Lambda}^{s}\left(  \beta,\omega\right)  $ commutes with
$E_{\Lambda}.$ The ground state of $\tilde{L}_{\Lambda}^{s}\left(
\beta,\omega\right)  ,$ that is to say
\[
\left\vert g_{\Lambda}\left(  \beta,\omega\right)  \right\rangle
:=\sum_{\alpha\subseteq\Lambda}g_{\alpha}^{\Lambda}\left(  \beta
,\omega\right)  \left\vert \alpha\right\rangle \,,
\]
where $g_{\alpha}^{\Lambda}\left(  \beta,\omega\right)  :=\frac{e^{-\frac
{\beta}{2}H_{\alpha}\left(  \omega\right)  }}{\left(  Z_{\Lambda}^{\left(
d\right)  }\left(  \beta,\omega\right)  \right)  ^{\frac{1}{2}}},$ belongs to
$\mathcal{H}_{\Lambda}^{+}.$
\end{remark}

Now let $\beta$ and $\omega$ be fixed. For any vector $\left\vert
u\right\rangle \in\mathcal{H}_{\Lambda},$ by (\ref{LbdsI}) the Dirichlet form
associated to $\tilde{L}_{\Lambda}^{s}\left(  \beta,\omega\right)  $ can be
written in the following way
\begin{align}
\left\langle u\right\vert \tilde{L}_{\Lambda}^{s}\left(  \beta,\omega\right)
\left\vert u\right\rangle  &  =\sum_{\alpha\subseteq\Lambda}\sum_{x\in\Lambda
}\frac{1}{\cosh\frac{\beta}{2}\Delta_{x}H_{\alpha}^{\omega}}\left[  \left(
\frac{u_{\alpha}-u_{\alpha\triangle\{x\}}}{2}\right)  ^{2}+u_{\alpha}^{2}%
\frac{e^{\frac{\beta}{2}\Delta_{x}H_{\alpha}^{\omega}}-1}{2}\right]
\label{formqL}\\
&  =\sum_{\alpha\subseteq\Lambda}\sum_{x\in\alpha}\frac{1}{\cosh\frac{\beta
}{2}\Delta_{x}H_{\alpha}^{\omega}}\left[  \left(  \frac{u_{\alpha}%
-u_{\alpha\backslash\{x\}}}{2}\right)  ^{2}+u_{\alpha}^{2}\frac{e^{\frac
{\beta}{2}\Delta_{x}H_{\alpha}^{\omega}}-1}{2}\right]  +\nonumber\\
& \quad +\sum_{\alpha\subseteq\Lambda}\sum_{x\in\alpha^{c}}\frac{1}{\cosh
\frac{\beta}{2}\Delta_{x}H_{\alpha}^{\omega}}\left[  \left(  \frac{u_{\alpha
}-u_{\alpha\cup\{x\}}}{2}\right)  ^{2}+u_{\alpha}^{2}\frac{e^{\frac{\beta}%
{2}\Delta_{x}H_{\alpha}^{\omega}}-1}{2}\right] \nonumber\\
&  =\sum_{\alpha\subseteq\Lambda}\sum_{x\in\alpha}\frac{1}{\cosh\frac{\beta
}{2}\Delta_{x}H_{\alpha}^{\omega}}\left[  \frac{\left(  u_{\alpha}%
-u_{\alpha\backslash\{x\}}\right)  ^{2}}{2}+u_{\alpha}^{2}\frac{e^{\frac
{\beta}{2}\Delta_{x}H_{\alpha}^{\omega}}-1}{2} \right. \nonumber \\
& \quad\left. +u_{\alpha\backslash\{x\}}%
^{2}\frac{e^{-\frac{\beta}{2}\Delta_{x}H_{\alpha}^{\omega}}-1}{2}\right]
\,.\nonumber
\end{align}
Clearly, $\forall\omega\in\Omega,\ \tilde{L}_{\Lambda}^{s}\left(
\beta=0,\omega\right)  =\bar{L}_{\Lambda}.$

\begin{proposition}
\label{prb}Let $\beta\geq0,\,\omega\in\Omega$ be fixed and $\Lambda$ be such
that $\left\vert \Lambda\right\vert =2N,\,N\in\mathbb{N}.$ For any $\left\vert
v\right\rangle \in\mathcal{H}_{\Lambda}^{\pm},$
\begin{equation}
\left\langle v\right\vert U_{\Lambda}\tilde{L}_{\Lambda}^{s}\left(
\beta,\omega\right)  U_{\Lambda}\left\vert v\right\rangle \leq\left(
1+2b_{dJ}\left(  \beta\right)  \right)  \left\langle v\right\vert L_{\Lambda
}\left\vert v\right\rangle \ , \label{rb}%
\end{equation}
where $b_{dJ}\left(  \beta\right)  $ is an analytic function of $\beta$ such
that $b_{dJ}\left(  \beta\right)  =c_{dJ}\beta+o\left(  \beta\right)  .$
\end{proposition}

\begin{proof}
Since by the previous remark it follows that $U_{\Lambda}\tilde{L}_{\Lambda
}^{s}\left(  \beta,\omega\right)  U_{\Lambda}$ commutes with $\bar{E}%
_{\Lambda},$ we can restrict ourselves to vectors in $\overline{\mathcal{H}%
}_{\Lambda}^{\pm}.$ Let us take $\left\vert v\right\rangle \in\overline
{\mathcal{H}}_{\Lambda}^{\pm},$\thinspace then \linebreak $U_{\Lambda}\left\vert
v\right\rangle =\left\vert u\right\rangle \in\mathcal{H}_{\Lambda}^{\pm}.$
From (\ref{formqL}) it follows that
\begin{align*}
\frac{1}{2}\sum_{\alpha\subseteq\Lambda}\sum_{x\in\alpha^{c}}\frac
{e^{\frac{\beta}{2}\Delta_{x}H_{\alpha}^{\omega}}-1}{\cosh\frac{\beta}%
{2}\Delta_{x}H_{\alpha}^{\omega}}u_{\alpha}^{2} & =\frac{1}{2}\sum_{\alpha
\subseteq\Lambda}\sum_{x\in\alpha^{c}}\frac{e^{\frac{\beta}{2}\Delta
_{x}H_{\alpha^{c}}^{\omega}}-1}{\cosh\frac{\beta}{2}\Delta_{x}H_{\alpha^{c}%
}^{\omega}}u_{\alpha^{c}}^{2}\\
& =\frac{1}{2}\sum_{\alpha^{c}\subseteq\Lambda}\sum_{x\in\alpha^{c}}%
\frac{e^{\frac{\beta}{2}\Delta_{x}H_{\alpha^{c}}^{\omega}}-1}{\cosh\frac
{\beta}{2}\Delta_{x}H_{\alpha^{c}}^{\omega}}u_{\alpha^{c}}^{2} \\
& =\frac{1}{2}
\sum_{\alpha\subseteq\Lambda}\sum_{x\in\alpha}\frac{e^{\frac{\beta}{2}%
\Delta_{x}H_{\alpha}^{\omega}}-1}{\cosh\frac{\beta}{2}\Delta_{x}H_{\alpha
}^{\omega}}u_{\alpha}^{2}\,.
\end{align*}
Thus,
\[
\left\langle u\right\vert \tilde{L}_{\Lambda}^{s}\left(  \beta,\omega\right)
\left\vert u\right\rangle =\sum_{\alpha\subseteq\Lambda}\sum_{x\in\alpha}%
\frac{1}{\cosh\frac{\beta}{2}\Delta_{x}H_{\alpha}^{\omega}}\left[
\frac{\left(  u_{\alpha}-u_{\alpha\backslash\{x\}}\right)  ^{2}}{2}+u_{\alpha
}^{2}\left(  e^{\frac{\beta}{2}\Delta_{x}H_{\alpha}^{\omega}}-1\right)
\right]  \,.
\]
Therefore,
\begin{align*}
\left\langle v\right\vert U_{\Lambda}\tilde{L}_{\Lambda}^{s}\left(
\beta,\omega\right)  U_{\Lambda}\left\vert v\right\rangle & =\left\langle
u\right\vert \tilde{L}_{\Lambda}^{s}\left(  \beta,\omega\right)  \left\vert
u\right\rangle \\
& \leq\left\langle v\right\vert L_{\Lambda}\left\vert
v\right\rangle +b_{J}\left(  \beta\right)  \left\langle u\right\vert
L_{\Lambda}\left\vert u\right\rangle \\
& =\left\langle v\right\vert L_{\Lambda}\left\vert v\right\rangle +b_{J}\left(
\beta\right)  \left\langle v\right\vert \bar{L}_{\Lambda}\left\vert
v\right\rangle \,,
\end{align*}
Moreover, if $v\in\mathcal{H}_{\Lambda}^{\pm}\cap\overline{\mathcal{H}%
}_{\Lambda}^{\pm},$ by (\ref{Lbar<L}),
\[
\left\langle v\right\vert L_{\Lambda}\left\vert v\right\rangle +b_{J}\left(
\beta\right)  \left\langle v\right\vert \bar{L}_{\Lambda}\left\vert
v\right\rangle \leq\left(  1+2b_{J}\left(  \beta\right)  \right)  \left\langle
v\right\vert L_{\Lambda}\left\vert v\right\rangle \,,
\]
where $b_{dJ}\left(  \beta\right)  :=\max_{z\in\lbrack0,4dJ]}\left[
\frac{e^{\frac{\beta}{2}z}-1}{\cosh\frac{\beta}{2}z}\right]  =\frac{e^{2\beta
dJ}-1}{\cosh2{\beta dJ}}.$
\end{proof}

In \cite{GI1, GI2} we introduced a new form for the generator of stochastic
Ising model with transition rates
\[
w^{\beta,\omega}\left(  \alpha,\alpha\triangle\left\{  x\right\}  \right)
=\frac{1+e^{\beta\Delta_{x}H_{\alpha}^{\omega}}}{4}%
\]
whose generic matrix element as an operator acting on $\mathcal{H}_{\Lambda}$
is

\begin{align}
& \left\langle \gamma\right\vert \hat{L}_{\Lambda}^{s}\left(  \beta
,\omega\right)  \left\vert \alpha\right\rangle  =\left\langle
\gamma\right\vert \sum_{x\in\Lambda}\hat{L}_{x}^{s}\left(  \beta
,\omega\right)  \left\vert \alpha\right\rangle \label{Lbds2}\\
&\left\langle \gamma\right\vert \hat{L}_{x}^{s}\left(  \beta,\omega\right)
\left\vert \alpha\right\rangle  \nonumber\\  =& \left\langle \gamma\right\vert \left\{
\mathbf{1}_{\alpha}\left(  x\right)  \cosh\frac{\beta}{2}\Delta_{x}H_{\alpha
}^{\omega}\left[  \bar{\ell}_{x}^{\Lambda}+I_{\Lambda}\frac{e^{\frac{\beta}%
{2}\Delta_{x}H_{\alpha}^{\omega}}-1}{2}\right]  +\right. \nonumber\\
& \qquad\quad \left.  +\mathbf{1}_{\alpha^{c}}\left(  x\right)  \cosh\frac{\beta}%
{2}\Delta_{x}H_{\alpha}^{\omega}\left[  \bar{\ell}_{x}^{\Lambda}+I_{\Lambda
}\frac{e^{\frac{\beta}{2}\Delta_{x}H_{\alpha}^{\omega}}-1}{2}\right]
\right\}  \left\vert \alpha\right\rangle \nonumber\\
&  =\left\langle \gamma\right\vert \left\{  \mathbf{1}_{\alpha}\left(
x\right)  \cosh\frac{\beta}{2}\Delta_{x}H_{\alpha}^{\omega}\left[
\frac{e^{\frac{\beta}{2}\Delta_{x}H_{\alpha}^{\omega}}+1}{2}\bar{\ell}%
_{x}^{\Lambda}+\frac{e^{\frac{\beta}{2}\Delta_{x}H_{\alpha}^{\omega}}-1}%
{2}\bar{\ell}_{x}^{\Lambda,\bot}\right]  +\right. \nonumber\\
& \qquad\quad \left.  +\mathbf{1}_{\alpha^{c}}\left(  x\right)  \cosh\frac{\beta}%
{2}\Delta_{x}H_{\alpha}^{\omega}\left[  \frac{e^{\frac{\beta}{2}\Delta
_{x}H_{\alpha}^{\omega}}+1}{2}\bar{\ell}_{x}^{\Lambda}+\frac{e^{\frac{\beta
}{2}\Delta_{x}H_{\alpha}^{\omega}}-1}{2}\bar{\ell}_{x}^{\Lambda,\bot}\right]
\right\}  \left\vert \alpha\right\rangle \,.\nonumber
\end{align}
We also showed that $\hat{L}_{\Lambda}^{s}\left(  \beta,\omega\right)  $
admits the representation
\[
\hat{L}_{\Lambda}^{s}\left(  \beta,\omega\right)  =\sum_{x\in\Lambda
}U_{\Lambda}e^{\frac{\beta}{2}\mathbf{H}_{\Lambda}\left(  \omega\right)  }%
\ell_{x}^{\Lambda}e^{-\beta\mathbf{H}_{\Lambda}\left(  \omega\right)  }%
\ell_{x}^{\Lambda}e^{\frac{\beta}{2}\mathbf{H}_{\Lambda}\left(  \omega\right)
}U_{\Lambda}\,,
\]
where
\begin{align}
\mathbf{H}_{\Lambda}\left(  \omega\right)   &  :=\sum_{\alpha\subseteq\Lambda
}H_{\alpha}\left(  \omega\right)  \left\vert \alpha\right\rangle \left\langle
\alpha\right\vert \simeq\mathbf{H}_{\Lambda}\left(  \omega\right)  =\sum
_{b\in\mathbb{B}_{\Lambda}}\omega_{b}\mathbf{s}_{b}\,,\label{defHop}\\
\mathbf{s}_{b}  &  :=\mathbf{1}_{b}\left(  x\right)  \mathbf{1}_{b}\left(
y\right)  \left(  1-\delta_{x,y}\right)  \mathbf{s}_{x}\mathbf{s}%
_{y}\,,\nonumber
\end{align}
so that, $\forall x\in\Lambda,\,\alpha\subseteq\Lambda,\,\mathbf{s}%
_{x}\left\vert \alpha\right\rangle =\left\vert \alpha\triangle
\{x\}\right\rangle $ (we prefer to work in the representation where the spin
flip operator is diagonal). By (\ref{LbdsI}) and (\ref{Lbds2}), for any two
basis vectors of $\mathcal{H}_{\Lambda},\ \left\vert \alpha\right\rangle
,\left\vert \gamma\right\rangle $ we have
\[
\left\vert \left\langle \gamma\right\vert \tilde{L}_{x}^{s}\left(
\beta,\omega\right)  \left\vert \alpha\right\rangle \right\vert \leq\left\vert
\left\langle \gamma\right\vert \hat{L}_{x}^{s}\left(  \beta,\omega\right)
\left\vert \alpha\right\rangle \right\vert \quad x\in\Lambda\,,
\]
moreover, the first order terms in the expansion for small $\beta$ of
$\left\langle \gamma\right\vert \tilde{L}_{\Lambda}^{s}\left(  \beta
,\omega\right)  \left\vert \alpha\right\rangle $ and $\left\langle
\gamma\right\vert \hat{L}_{\Lambda}^{s}\left(  \beta,\omega\right)  \left\vert
\alpha\right\rangle $ are equal for every $\Lambda.$ Clearly, $\hat
{L}_{\Lambda}^{s}\left(  \beta,\omega\right)  $ also commutes with
$E_{\Lambda}$ and for any $\left\vert u\right\rangle \in\mathcal{H}_{\Lambda
},$ since $\left\vert u\right\rangle =\left\vert u^{+}\right\rangle
+\left\vert u^{-}\right\rangle ,\,\left\vert u^{\pm}\right\rangle
\in\mathcal{H}_{\Lambda}^{\pm},$ we have
\begin{align}
\left\langle u^{\pm}\right\vert \hat{L}_{\Lambda}^{s}\left(  \beta
,\omega\right)  \left\vert u^{\pm}\right\rangle  &  =\sum_{\alpha
\subseteq\Lambda}\sum_{x\in\Lambda}\cosh\frac{\beta}{2}\Delta_{x}H_{\alpha
}^{\omega}\left[  \left(  \frac{u_{\alpha}^{\pm}-u_{\alpha\triangle\{x\}}%
^{\pm}}{2}\right)  ^{2}+\left(  u_{\alpha}^{\pm}\right)  ^{2}\frac
{e^{\frac{\beta}{2}\Delta_{x}H_{\alpha}^{\omega}}-1}{2}\right]  \label{formL2}%
\\
&  =\sum_{\alpha\subseteq\Lambda}\sum_{x\in\alpha}\cosh\frac{\beta}{2}%
\Delta_{x}H_{\alpha}^{\omega}\left[  \left(  \frac{u_{\alpha}^{\pm}%
-u_{\alpha\backslash\{x\}}^{\pm}}{2}\right)  ^{2}+\left(  u_{\alpha}^{\pm
}\right)  ^{2}\frac{e^{\frac{\beta}{2}\Delta_{x}H_{\alpha}^{\omega}}-1}%
{2}\right]  \,.\nonumber
\end{align}
Proceding as in Proposition \ref{prb},\ we get
\begin{equation}
\left\langle u\right\vert \hat{L}_{\Lambda}^{s}\left(  \beta,\omega\right)
\left\vert u\right\rangle \leq\left(  1+2b_{dJ}^{\prime}\left(  \beta\right)
\right)  \left\langle u\right\vert \bar{L}_{\Lambda}\left\vert u\right\rangle
\qquad u\in\mathcal{H}_{\Lambda}^{\pm}\cap\overline{\mathcal{H}}_{\Lambda
}^{\pm}\,, \label{rb2}%
\end{equation}
where $b_{dJ}^{\prime}\left(  \beta\right)  :=\max_{z\in\lbrack0,4dJ]}\left[
\left(  e^{\frac{\beta}{2}z}-1\right)  \cosh\frac{\beta}{2}z+\cosh\frac{\beta
}{2}z-1\right]  =\frac{e^{2\beta dJ}-1}{2}.$ Comparing the Dirichlet forms of
$\hat{L}_{\Lambda}^{s}\left(  \beta,\omega\right)  $ and $\tilde{L}_{\Lambda
}^{s}\left(  \beta,\omega\right)  ,$ we proved in \cite{GI2} that this process
converges to the equilibrium state at high temperature faster than the
heat-bath process.

\begin{remark}
\label{r1}The relative bounds (\ref{rb}) and (\ref{rb2}) are independent of
$\Lambda$ and extend straightforwardly to the quadratic forms associated to
the operators $\tilde{L}^{s}\left(  \beta,\omega\right)  $ and $\hat{L}%
^{s}\left(  \beta,\omega\right)  $ acting on $\mathcal{H}.$ Therefore, by
standard arguments of perturbation theory (see \cite{K} Theorem VI.3.4)
(\ref{rb}) implies the analyticity of the projectors
\[
P_{n}\left(  \beta,\omega\right)  :=\oint_{\{z\in\mathbb{C}\,:\,\left\vert
z-n\right\vert \leq r\left(  \beta\right)  \}}\frac{dz}{2\pi i}\frac
{1}{Iz-A\left(  \beta,\omega\right)  }\quad n\in\mathbb{N}\ ,
\]
where $A\left(  \beta,\omega\right)  $ is either $\tilde{L}^{s}\left(
\beta,\omega\right)  $ or $\hat{L}^{s}\left(  \beta,\omega\right)  ,$ for
sufficiently small values of $\beta.$
\end{remark}

\subsubsection{Lower bound $g_{d}^{-}\left(  \beta\right)  $}

By Remark 3, we can make use of perturbation theory and, for sufficently small
values of $\beta$ and any realization of the potential, we can write
\[
\hat{L}_{\Lambda}^{s}\left(  \beta,\omega\right)  =U_{\Lambda}L_{\Lambda
}U_{\Lambda}+\beta U_{\Lambda}T_{\Lambda}^{\left(  1\right)  }\left(
\omega\right)  U_{\Lambda}+\bar{T}_{\Lambda}\left(  \beta,\omega\right)  \,,
\]
where
\[
T_{\Lambda}^{\left(  1\right)  }\left(  \omega\right)  :=\frac{1}{2}\sum
_{x\in\Lambda}[[\mathbf{H}_{\Lambda}\left(  \omega\right)  ,\ell_{x}^{\Lambda
}],\ell_{x}^{\Lambda}]
\]
is the first term in the expansion of $\hat{L}_{\Lambda}^{s}\left(
\beta,\omega\right)  $ and $\bar{T}_{\Lambda}\left(  \beta,\omega\right)  $ is
such that
\[
\left\langle u\right|  \bar{T}_{\Lambda}\left(  \beta,\omega\right)  \left|
u\right\rangle \leq\beta^{2}C\left(  d,J\right)  \,,
\]
with $C\left(  d,J\right)  $ a positive constant.

Since, by definition of $U_{\Lambda}$, $U_{\Lambda}\hat{L}_{\Lambda}%
^{s}\left(  \beta,\omega\right)  U_{\Lambda}$ and $\hat{L}_{\Lambda}%
^{s}\left(  \beta,\omega\right)  $ have the same spectrum, the eigenspace
corresponding to the first non-trivial eigenvalue of the unperturbed
generator, $\xi_{1}\left(  L_{\Lambda}\right)  =1,$ is $span\{\left\vert
y\right\rangle :y\in\Lambda\}$ and
\begin{align}
\left\langle z\right\vert T_{\Lambda}^{\left(  1\right)  }\left(
\omega\right)  \left\vert y\right\rangle  &  =\frac{1}{2}\left\langle
z\right\vert \sum_{x\in\Lambda}[[\mathbf{H}_{\Lambda}\left(  \omega\right)
,\ell_{x}^{\Lambda}],\ell_{x}^{\Lambda}]\left\vert y\right\rangle
\label{pertT}\\
&  =\frac{1}{2}\sum_{x\in\Lambda}\left\langle z\right\vert \mathbf{H}%
_{\Lambda}\left(  \omega\right)  \left\vert y\right\rangle \left(
\delta_{x,y}+\delta_{x,z}-2\delta_{z,x}\delta_{x,y}\right) \nonumber\\
&  =\left\langle z\right\vert \mathbf{H}_{\Lambda}\left(  \omega\right)
\left\vert y\right\rangle -\delta_{z,y}\left\langle y\right\vert
\mathbf{H}_{\Lambda}\left(  \omega\right)  \left\vert y\right\rangle
\,,\nonumber
\end{align}
where by (\ref{defHop})
\begin{equation}
\left\langle z\right\vert \mathbf{H}_{\Lambda}\left(  \omega\right)
\left\vert y\right\rangle =\sum_{b\in\mathbb{B}_{\Lambda}}\omega
_{b}\left\langle z\right\vert \mathbf{s}_{b}\left\vert y\right\rangle
=\sum_{b\in\mathbb{B}_{\Lambda}}\omega_{b}\left\langle z|\{y\}\triangle
b\right\rangle =\sum_{b\in\mathbb{B}_{\Lambda}}\omega_{b}\mathbf{1}%
_{\{z,y\}}\left(  b\right)  =\omega_{z,y}\,. \label{pertH}%
\end{equation}

Moreover, looking at the expansion in $\beta$ of the Dirichlet forms of
$\tilde{L}_{\Lambda}^{s}\left(  \beta,\omega\right)  $ and $\hat{L}_{\Lambda
}^{s}\left(  \beta,\omega\right)  ,$ we realize that these operators coincides
up to first order. Hence, we get
\[
\xi_{1}\left(  \tilde{L}_{\Lambda}^{s}\left(  \beta,\omega\right)  \right)
\geq g_{d}^{-}\left(  \beta\right)  \,,
\]
with $g_{d}^{-}\left(  \beta\right)  $ analytic function of $\beta$ such that
\begin{equation}
g_{d}^{-}\left(  \beta\right)  :=1-\beta\sup_{z\in\Lambda}\sum_{y\in\Lambda
}\left\vert \omega_{x,y}\right\vert +o\left(  \beta\right)  =1-2\beta
dJ+O\left(  \beta^{2}\right)  \,. \label{lb}%
\end{equation}
Notice that all the above estimates, which are independent of $\Lambda$, hold
in infinite volume as well.

\begin{remark}
Since the $\omega$'s are bounded, the last result implies the existence of a
value of $\beta_{d}\left(  J\right)  $ smaller than the critical one
$\beta_{c}\left(  d,\omega\right)  ,$ such that for $\mathbb{P}$ a.e.
$\omega,$ if $\beta\in\lbrack0,\beta_{d}\left(  J\right)  ),$ the process is
ergodic. Hence, by the reversibility with respect to the Gibbs measure, we get
the uniqueness of the Gibbs state. Furthermore, the unique element $\mu
^{\beta,\omega}$ of $\mathcal{G}\left(  \beta,\omega\right)  $ has\thinspace
the property
\begin{equation}
\mu^{\beta,\omega}\left(  A\right)  =\mu^{\beta,\theta_{z}\omega}\left(
\tau_{z}A\right)  \quad A\subset\mathcal{S},\,z\in\mathbb{Z}^{d}\,,
\label{invtr}%
\end{equation}
where
\[
\tau_{z}A:=\left\{  \sigma\in\mathcal{S}:\forall x\in\mathbb{Z}^{d}\quad
\sigma_{x}=\eta_{x-z}=\left(  \tau_{z}\eta\right)  _{x},\,\eta\in A\right\}
\,.
\]

\end{remark}

Let $\left\{  \Theta_{z}\right\}  _{z\in\mathbb{Z}^{d}}$ be the unitary group
of operators on $\mathcal{L}\left(  \beta,\omega\right)  $ generated by the
group $\left\{  \tau_{z}\right\}  _{z\in\mathbb{Z}^{d}}$ that is,
\[
\left(  \Theta_{z}\varphi\right)  \left(  \sigma\right)  =\varphi\left(
\tau_{z}^{-1}\sigma\right)  \quad\varphi\in\mathcal{L}\left(  \beta
,\omega\right)  \,.
\]
Then,\ by the previous remark, we get that $\forall z\in\mathbb{Z}^{d}$ the
Hilbert spaces $\mathcal{L}\left(  \beta,\omega\right)  $ and $\mathcal{L}%
\left(  \beta,\theta_{z}^{-1}\omega\right)  $ are unitary equivalent
(isomorphic) via the unitary mapping $\Theta_{z}$
\[
\Theta_{z}:\mathcal{L}\left(  \beta,\omega\right)  \longmapsto\mathcal{L}%
\left(  \beta,\theta_{z}^{-1}\omega\right)
\]
and from the representation (\ref{defgen}) of $L\left(  \beta,\omega\right)
,$ we have
\[
\Theta_{z}L\left(  \beta,\omega\right)  \Theta_{z}^{-1}=L\left(  \beta
,\theta_{z}^{-1}\omega\right)  \,,
\]
which implies that, at least for $\beta\in\left[  0,\beta_{d}\left(  J\right)
\right)  ,$ the family of operators and spaces $\left(  L\left(  \beta
,\omega\right)  ,\mathcal{L}\left(  \beta,\omega\right)  \right)  $ is a
metrically transitive family with respect to the unitary group of lattice
translation $\left\{  \Theta_{z}\right\}  _{z\in\mathbb{Z}^{d}}.$ Hence, (see
\cite{PF} and Remark 4 of \cite{AMSZ}) the spectrum of $L\left(  \beta
,\omega\right)  $ is non-random for $\mathbb{P}$-a.e. $\omega.$

\begin{remark}
To get an upper bound for the spectral gap of the generator of the process we
can compute the Dirichlet form of $L\left(  \beta,\omega\right)  $ with
respect to the function of the empirical magnetization
\[
\phi_{\Lambda}:=\sum_{x\in\Lambda}\frac{\sigma_{x}}{\left|  \Lambda\right|
}-\mu^{\beta,\omega}\left(  \sum_{x\in\Lambda}\frac{\sigma_{x}}{\left|
\Lambda\right|  }\right)  \,.
\]

We have
\begin{align*}
\left\langle \phi_{\Lambda},L\left(  \beta,\omega\right)  \phi_{\Lambda
}\right\rangle _{\beta,\omega}  &  =\frac{1}{2}\int\mu^{\beta,\omega}\left(
d\sigma\right)  \sum_{x\in\mathbb{Z}^{d}}w_{x}^{\beta,\omega}\left(
\sigma\right)  \left[  \sum_{y\in\Lambda}\frac{\sigma_{y}}{\left\vert
\Lambda\right\vert }\left(  1-2\delta_{x,y}\right)  -\sum_{y\in\Lambda}%
\frac{\sigma_{y}}{\left\vert \Lambda\right\vert }\right]  ^{2}\\
&  =2\int\mu^{\beta,\omega}\left(  d\sigma\right)  \sum_{x\in\Lambda}%
\frac{w_{x}^{\beta,\omega}\left(  \sigma\right)  }{\left\vert \Lambda
\right\vert ^{2}}\,.
\end{align*}
Dividing by the $\mathcal{L}\left(  \beta,\omega\right)  $ norm of
$\phi_{\Lambda}$
\[
\frac{1}{\left\vert \Lambda\right\vert ^{2}}\sum_{x,y\in\Lambda}\left[
\mu^{\beta,\omega}\left(  \sigma_{x}\sigma_{y}\right)  -\mu^{\beta,\omega
}\left(  \sigma_{x}\right)  \mu^{\beta,\omega}\left(  \sigma_{y}\right)
\right]  \ ,
\]
we have that the spectral gap is smaller than
\[
\frac{2\int\mu^{\beta,\omega}\left(  d\sigma\right)  \sum_{x\in\Lambda}%
w_{x}^{\beta,\omega}\left(  \sigma\right)  }{\sum_{x,y\in\Lambda}\left[
\mu^{\beta,\omega}\left(  \sigma_{x}\sigma_{y}\right)  -\mu^{\beta,\omega
}\left(  \sigma_{x}\right)  \mu^{\beta,\omega}\left(  \sigma_{y}\right)
\right]  }\,.
\]
By the ergodicity of the random field $\omega$ with respect to the lattice
translations, the last expression becomes
\begin{equation}
\frac{2\int\mathbb{P}\left(  d\omega\right)  \int\mu^{\beta,\omega}\left(
d\sigma\right)  w_{0}^{\beta,\omega}\left(  \sigma\right)  }{\int
\mathbb{P}\left(  d\omega\right)  \sum_{y\in\mathbb{Z}^{d}}\left[  \mu
^{\beta,\omega}\left(  \sigma_{0}\sigma_{y}\right)  -\mu^{\beta,\omega}\left(
\sigma_{0}\right)  \mu^{\beta,\omega}\left(  \sigma_{y}\right)  \right]  }\,,
\label{ubgap}%
\end{equation}
where
\[
\chi^{d,\omega}\left(  \beta\right)  :=\sum_{x\in\mathbb{Z}^{d}}\left[
\mu^{\beta,\omega}\left(  \sigma_{x}\sigma_{0}\right)  -\mu^{\beta,\omega
}\left(  \sigma_{x}\right)  \mu^{\beta,\omega}\left(  \sigma_{0}\right)
\right]  \ ,\quad\omega\in\Omega\ ,
\]
is the \emph{susceptibility} relative to a realization of the potential. We
could now get estimates for (\ref{ubgap}) at small values of $\beta$ through a
cluster expansion. We will not pursue this here but rather get a bound by
different means in 3.1.2 below. In the ferromagnetic case $(\omega_{b}\geq
J^{-}>0,\forall b\in\mathbb{B}_{d}),$ by the Griffiths inequalities
(see for example \cite{L} page 186), we have that $\chi^{d,\omega}\left(
\beta\right)  $ is larger than or equal to the susceptibility relative to the
configuration of the potential constantly equal to $J^{-},$ $\chi^{d,J^{-}%
}\left(  \beta\right)  ,$ which is known to be a function of $\beta$ diverging
when $\beta$ approachs its critical value $\beta_{c}\left(  d,\omega\right)  $
from below. In particular, in the two-dimensional case, $\chi^{2,J^{-}}\left(
\beta\right)  $ is proportional to $\left\vert \beta-\beta_{c}\left(
2,J^{-}\right)  \right\vert ^{-\frac{7}{4}}$ ( \cite{H} Theorem 2.11). Then by
(\ref{stimw}), (\ref{ubgap}) is smaller than
\[
\frac{2\psi\left(  -\beta4dJ\right)  \vee\psi\left(  \beta4dJ\right)  }%
{\chi^{d,J^{-}}\left(  \beta\right)  }\,.
\]

\end{remark}

\subsubsection{Upper bound $g_{d}^{+}\left(  \beta\right)  $}

Since, for any $\beta\geq0$ and $\omega\in\Omega,$
\[
\left\langle u\right\vert \tilde{L}_{\Lambda}^{s}\left(  \beta,\omega\right)
\left\vert u\right\rangle \leq\left\langle u\right\vert \hat{L}_{\Lambda}%
^{s}\left(  \beta,\omega\right)  \left\vert u\right\rangle, \quad\left\vert
u\right\rangle \in\mathcal{H}_{\Lambda}\,,
\]
from (\ref{pertT}) and (\ref{pertH}) we get $\xi_{1}\left(  \tilde{L}%
_{\Lambda}^{s}\left(  \beta,\omega\right)  \right)  \leq g_{d}^{+}\left(
\beta\right)  ,$ with $g_{d}^{+}\left(  \beta\right)  $ analytic function of
$\beta$ such that
\begin{equation}
g_{d}^{+}\left(  \beta\right)  :=1+\beta\sup_{z\in\Lambda}\sum_{y\in\Lambda
}\left\vert \omega_{x,y}\right\vert +o\left(  \beta\right)  =1+2\beta
dJ+O\left(  \beta^{2}\right)  \,. \label{ub}%
\end{equation}

\subsubsection{$\sigma_{\beta}^{\left(  1\right)  }\subseteq\left[  g_{d}%
^{-}\left(  \beta J\right)  ,g_{d}^{+}\left(  \beta J\right)  \right]  $}

We just notice that, for $\beta$ smaller than $\beta_{d}\left(  J\right)  ,$
we have $g_{d}^{-}\left(  \beta\right)  =1-2\beta Jd+O\left(  \beta
^{2}\right)  $ and $g_{d}^{+}\left(  \beta\right)  =1+2\beta Jd+O\left(
\beta^{2}\right)  ,$ which implies that, for such values of $\beta,$
\begin{equation}
\sigma_{\beta}^{\left(  1\right)  }\subseteq\left[  1-2\beta Jd,1+2\beta
Jd\right]  \,. \label{>}%
\end{equation}

\subsection{Proof of Theorem \ref{mt}}

Here we mimic the second part of the proof of Theorem 3 in \cite{AMSZ} and
consider $\Omega$ as a topological space endowed with the Schwartz topology,
which we will denote by $\mathcal{D}_{\mathbb{B}_{d}}.$ We will denote by
$supp\mathbb{P}$ the support of $\mathbb{P}$ as a function on $\mathcal{D}%
_{\mathbb{B}_{d}}.$ Let $\zeta$ be any realization of the potential constant
on $\mathbb{B}_{d}$ which belongs to the support of $\mathbb{P},$ namely
\[
\zeta=\left\{  \omega_{b}=\zeta\in\mathbb{R},\,\forall b\in\mathbb{B}\right\}
\in supp\mathbb{P}%
\]
and denote by $\mathcal{C}_{\mathbb{B}_{d}}\subset\mathcal{D}_{\mathbb{B}_{d}%
}$ the collection of all such realizations of the potential.

Theorem 3 of \cite{AMSZ} uses the explicit representation of the matrix
elements of the generator for the one dimension model to prove the weak
continuity of the spectral measure. All is really needed is that the matrix
elements of the generator and thus the semigroup are smooth functions of the
potential. In higher dimension we rely on (\ref{LbdsI}), which in particular
ensures the necessary regularity.

It is proved in Theorem 2.2 of \cite{M} that, for any constant realizations
$\zeta$ of the potential in $supp\mathbb{P},$ there exists a value $\beta
_{d}^{\left(  1\right)  }\left(  \zeta\right)  >0$ such that, for any
$\left\vert \beta\right\vert <\beta_{d}^{\left(  1\right)  }\left(
\zeta\right)  $, we obtain
\[
\left[  1-a_{d}\left(  \beta\zeta\right)  ,1+a_{d}\left(  \beta\zeta\right)
\right]  \subseteq\sigma_{\beta}^{\left(  1\right)  }\,,
\]
with
\begin{equation}
a_{d}\left(  r\right)  :=\max_{\lambda\in\mathbb{T}^{d}}\left\vert
a_{d}\left(  \lambda,r\right)  \right\vert \,. \label{adbz}%
\end{equation}
For the definition of $a_{d}\left(  \lambda,r\right)  $ see Theorem 2.3 of
\cite{M}. Consequently, if \newline$\beta_{d}^{\left(  1\right)  }:=\inf_{\zeta}%
\beta_{d}^{\left(  1\right)  }\left(  \zeta\right)  >0$ and $\beta_{d}^{\ast
}\left(  J\right)  :=\beta_{d}\left(  J\right)  \wedge\beta_{d}^{\left(
1\right)  },$ then, $\forall\beta\in\left[  0,\beta_{d}^{\ast}\left(
J\right)  \right)  ,$
\begin{equation}
\left[  1-\bar{a}_{d}\left(  \beta\right)  ,1+\bar{a}_{d}\left(  \beta\right)
\right]  \subseteq\bigcup_{\zeta\in\mathcal{C}_{\mathbb{B}_{d}}}\sigma_{\beta
}^{\left(  1\right)  }\left(  L^{\left(  1\right)  }\left(  \beta
,\zeta\right)  \right)  \subseteq\sigma_{\beta}^{\left(  1\right)  }\ ,
\label{<}%
\end{equation}
with
\[
\bar{a}_{d}\left(  \beta\right)  :=\max_{\zeta\in\mathcal{C}_{\mathbb{B}_{d}}%
}a_{d}\left(  \beta\zeta\right)  \ ,
\]
which is an analytic function of $\beta.$ Thus $f_{d}^{\pm}\left(
\beta\right)  :=\pm\bar{a}_{d}\left(  \beta\right)  .$ Since $\bar{a}%
_{d}(\beta)=a_{d}(\beta\zeta)$ for $\zeta$ such that $\left\vert
\zeta\right\vert =J,$ then for small values of $\beta,$ $\bar{a}_{d}\left(
\beta\right)  =2dJ\beta+o\left(  \beta\right)  ,$ where the linear term in
$\beta$ is the same in the expansion of (\ref{lb}), as well as in (\ref{ub}).

\subsection{Proof of Theorem \ref{ct}}

Because the family of operators and spaces \newline$\left(  L\left(  \beta
,\omega\right)  ,\mathcal{L}\left(  \beta,\omega\right)  \right)  $ is
metrically transitive with respect to lattice translations, $\sigma_{\beta
}^{\left(  1\right)  }$ is a non random set (see Remark 2). Thus, for every
$\beta\in\left[  0,\beta_{d}^{\ast}\left(  J\right)  \right)  ,$ at first
order in $\beta,$ by (\ref{>}) and (\ref{<}) we obtain
\[
\left[  1-2dJ\beta,1+2dJ\beta\right]  \subseteq\sigma_{\beta}^{\left(
1\right)  }\subseteq\left[  1-2\beta dJ,1+2\beta dJ\right]  \,.
\]

\ 

\noindent\textbf{Acknowledgements.} The authors thank R. Minlos for his
encouragement and an anonymous referee for useful comments and for pointing
out an incorrect formulation of Theorem \ref{ct} in an earlier draft and
other useful comments which led to an improvement of the paper.

\end{document}